\documentclass[english,12pt,leqno]{article}
\usepackage{amsxtra}
\usepackage{amsmath}
\usepackage{amssymb}
\usepackage{amsfonts}
\usepackage[all]{xy}
\usepackage{mathrsfs}
\usepackage{amsthm}
\usepackage{enumitem}
\usepackage{mathabx}
\usepackage{bbm}
\usepackage{dsfont}
\usepackage{esint}    
\usepackage{subcaption}
\usepackage{booktabs}
\usepackage{tabularx}
\usepackage{rotating}
\usepackage{appendix}
\usepackage{graphicx, nicefrac}
\let\rfb\reflectbox

\newcommand{\uglyfrac}[2]{\rfb{\nicefrac{\rfb{#1}}{\rfb{#2}}}}

\newcommand{\muglyfrac}[2]{\uglyfrac{\ensuremath{#1}}{\ensuremath{#2}}}

\usepackage{color}
\usepackage{textcomp}
\usepackage{geometry}

\usepackage{tikz-cd}
\usepackage{algorithm} 
\usepackage[noend]{algpseudocode} 
\usepackage[pdf]{pstricks}

\psset{unit=1pt}
\psset{arrowsize=4pt 1}
\psset{linewidth=.5pt}

\usepackage{babel}

\usepackage[final]{hyperref}

\makeatletter
\hypersetup{
    unicode=true,           
    pdftoolbar=true,        
    pdfmenubar=true,        
    pdffitwindow=false,     
    pdfstartview={FitH},    
    pdftitle={\@title},     
    pdfauthor={\@author},   
    pdfsubject={},          
    pdfcreator={},          
    pdfproducer={},         
    pdfkeywords={},         
    pdfnewwindow=true,      
    colorlinks,             
    linkcolor=black,        
    citecolor=black,        
    filecolor=black,        
    urlcolor=black          
}
\makeatother
 
\newtheorem{thm}{Theorem}[section]
\newtheorem{cor}[thm]{Corollary}
\newtheorem{lem}[thm]{Lemma}
\newtheorem{prop}[thm]{Proposition}

\newtheoremstyle{example}{\topsep}{\topsep}%
     {}
     {}
     {\bfseries}
     {.}
     {2pt}
     {\thmname{#1}\thmnumber{ #2}\thmnote{ #3}}

   \theoremstyle{example}
   
   \newtheorem{Defi}[thm]{Definition}
   \newtheorem{rem}[thm]{Remark}

   \newtheorem{ex}[thm]{Example}

\numberwithin{equation}{thm}
\DeclareMathOperator{\Spec}{Spec}
\DeclareMathOperator{\PExp}{PExp}

\DeclareMathOperator{\GW}{GW}
\DeclareMathOperator{\Hom}{Hom}
\DeclareMathOperator{\MW}{MW}
\DeclareMathOperator{\nicedeg}{deg}
\DeclareMathOperator{\chern}{ch}
\DeclareMathOperator{\Trop}{Trop}
\DeclareMathOperator{\nicespan}{span}
\DeclareMathOperator{\Rel}{Rel}
\DeclareMathOperator{\Tor}{Tor}
\DeclareMathOperator{\relint}{relint}
\DeclareMathOperator{\forgetful}{forgetful}

\def\CC{\mathbb{C}}

\def\PP{\mathbb{P}}
\def\RR{\mathbb{R}}
\def\ZZ{\mathbb{Z}}
\def\QQ{\mathbb{Q}}

\renewcommand{\precdot}{\prec\mathrel{\mkern-5mu}\mathrel{\cdot}}

\newcommand\scalemath[2]{\scalebox{#1}{\mbox{\ensuremath{\displaystyle #2}}}} 

\newcommand{\opk}{\mathrm{op}K} 
\newcommand{\Rho}{\mathrm{P}}

\def\Oc{\mathcal{O}}

 \def\codim{{\on{codim}}}
 \def\dim{{\on{dim}}}

 \def\Ext{{\on{Ext}}}

\usepackage{calligra}

  \def\Hom{\operatorname{Hom}\nolimits}

\def\mult{\on{mult}}

\def\on{\operatorname}

 \def\Spec {\on{Spec}}

\def\wt{\widetilde}

\def\1{{\mathbf{1}}}

\def\(({(\hskip -1mm (}
\def\)){)\hskip -1mm )}
 \def\be{\begin{equation}}
\def\ee{\end{equation}}
\def\ed{\end{document}}

\title{Minkowski weights and the Grothendieck group of a toric variety}
\author{Aniket Shah}
 
\begin{document}

\maketitle

\begin{abstract} For a fan $\Delta$, we introduce Grothendieck weights as a ring of functions from $\Delta$ to $\ZZ$ that form a K-theoretic analogue of Minkowski weights and describe the operational $K$-theory of a complete toric variety. We give an explicit balancing condition and product formula for these weights, and describe relationships with other fan-based invariants. Applications are given to vector bundles on a toric surface, and to the calculation of Euler characteristics on sch\" on subvarieties.
\end{abstract}
 
\addtocounter{section}{1}

\section*{Introduction}
The operational Chow ring of an algebraic variety $X$ is a way of extending the classical Chow ring of cycles modulo rational equivalence to singular varieties. From here on we restrict to the case when $X$ is a toric variety over $\CC$, in which case $X$ is equivalent to the data of a rational polyhedral fan $\Delta$. In \cite{FS}, Fulton and Sturmfels showed that when $X$ is complete the operational Chow ring $A^*(X)$ is isomorphic to a ring of balanced $\ZZ$-valued functions on $\Delta$, which they called \textit{Minkowski weights}, with a ``displacement rule" for calculating products. Several years later, the equivariant theory was described by Payne \cite{P}, and it was shown by Katz and Payne in \cite{KP} that the equivariant Chow ring $A^*_T(X)$ does not surject onto $A^*(X)$ in general. 

We consider the operational $K$-theory ring $\opk^\circ(X)$ defined by Anderson and Payne in \cite{AP}. This ring is isomorphic to the Grothendieck ring of vector bundles $K^\circ(X)$ when $X$ is smooth, but may be thought of as a better-behaved alternative when $X$ is singular - to see how $K^\circ(X)$ may be pathological even when $X$ is toric, see \cite{G}. Anderson and Payne gave a description of the equivariant theory $\opk^\circ_T(X)$ in terms of piecewise exponential functions on $\Delta$, but just as for Chow groups the forgetful map $\opk^\circ_T(X)\rightarrow\opk^\circ(X)$ is not always surjective: a $3$-fold example is given in \cite[Section 7]{AGP} over $\QQ$, and in this paper we give a surface example over $\ZZ$ (see Corollary \ref{cor:nonsurj}).

Since $A^*(X),A^*_T(X)$, and $\opk^\circ_T(X)$ have combinatorial descriptions, it remains to describe $\opk^\circ(X)$. We show that $\opk^\circ(X)$ is isomorphic to another ring of $\ZZ$-valued functions on $\Delta$, which we call \textit{Grothendieck weights}. These functions are also balanced in a sense that can be stated via Ehrhart theory, e.g. Proposition \ref{prop:ehrhartequiv}. Our primary goal here is to study properties of these weights.

Our secondary goal is to connect Grothendieck weights to tropical geometry. There, one starts with a valued field $K$ and a subvariety $Y$ inside $(K^*)^n$, and obtains a subset of $\RR^n$ which is the support of a rational polyhedral complex $\Trop(Y)$. Let us fix $K=\CC$ with the trivial valuation for simplicity, in which case $\Trop(Y)$ is simply the support of a rational polyhedral fan. By Tevelev's theorem \cite{Te}, the closure of $Y$ in a toric variety $X(\Delta)$ partially compactifying $(\CC^*)^n$ is proper precisely when the fan $\Delta$ contains $\Trop(Y)$. If Y is \textit{sch\" on}, then one may additionally define a function $f_{\overline{Y}}$ from the cones of $\Delta$ to $\ZZ$ which remembers intersection-theoretic properties of $\overline{Y}$ (see e.g. \cite{Mi,AR,K}), and is a Minkowski weight. The methods of tropical intersection theory explain how one can use the function $f_{\overline{Y}}$ to compute intersection numbers of divisors on $Y$. This has been carried out in many contexts, for example by Kerber and Markwig in \cite{KM} to compute intersection numbers of Psi classes on $\overline{M}_{0,n}$. Here, we will define a function $g_{\overline{Y}}$ from $\Delta$ to $\ZZ$ which remembers $K$-theoretic properties of $\overline{Y}$. This function will be balanced, but in the sense of Grothendieck weights rather than Minkowski weights.

Since $X$ is toric, the dual of $K_\circ(X)$ is a ring, by Proposition \ref{prop:isring}. We develop the ring of Grothendieck weights in Section $3$ to describe this ring. For Minkowski weights, the balancing condition comes from relations between fundamental classes of $T$-invariant subvarieties in $A_*(X)$, and analogously the balancing condition for Grothendieck weights will come from relations between the classes of structure sheaves of $T$-invariant subvarieties $K_\circ(X)$. However, determining these relations directly is equivalent to a problem in Ehrhart theory which seems fairly difficult, see Proposition \ref{prop:ehrhartequiv}. 

When $X$ is complete and toric, the operational K-theory ring $\opk^\circ(X)$ introduced in \cite{AP} coincides with the dual of $K_\circ(X)$, see \cite[Theorem 1.3]{AP}. Thus, simply due to general properties of operational K-theory, one obtains maps between Grothendieck weights and other well-known invariants of toric fans, namely Minkowski weights and piecewise exponential functions on $\Delta$. 

Our main results on Grothendieck weights are: 
\begin{itemize}
	\item a conceptual characterization of Grothendieck weights on projective simplicial fans in terms of Ehrhart polynomials (Proposition \ref{prop:ehrhartequiv});
	\item an explicit ``balancing condition" characterization of Grothendieck weights, based on a Riemann-Roch matrix we compute using previous work of \cite{PT} (Theorem \ref{thm:balancing});
	\item a displacement rule for products of Grothendieck weights (Theorem \ref{thm:prod});
	\item expressions for the maps to Grothendieck weights in the complete case (Section \ref{sec:mapstogw});
	\item an explicit computation demonstrating that there exists a complete toric surface such that the natural map $K^\circ_T(X)\rightarrow K^\circ(X)$ is not surjective (Corollary \ref{cor:nonsurj}).
\end{itemize}
With these properties of Grothendieck weights established we return in the last section to the tropical context to demonstrate that it is possible to use $g_{\overline{Y}}$ and our explicit formulae for these maps to compute classical Euler characteristics of line bundles or more generally coherent sheaves on a sch\" on subvariety $Y$. To that end, we describe an example computation on $\overline{M}_{0,n}$.

\subsection*{Notation and Conventions}

We establish our notation and conventions regarding toric varieties--generally we follow \cite{F93}. By toric variety we mean a normal variety with a $(\CC^*)^n$-action, such that there is a dense orbit isomorphic to $(\CC^*)^n$. We write $T$ for $(\CC^*)^n$, and write $M$ for the character lattice $\Hom_{alg. gp.}(T,\CC^*)$ and $N\cong M^\vee$ for the lattice of one parameter subgroups $\Hom_{alg. gp.}(\CC^*,T)$. 

As is well-known, $X$ corresponds to a polyhedral fan $\Delta$, which is a finite set of strongly convex rational polyhedral cones in $N$, such that for $\alpha\in\Delta$, any face of $\alpha$ is also in $\Delta$, and for any two cones $\alpha$ and $\beta$ in $\Delta$, the intersection $\alpha\cap\beta$ is a face of both $\alpha$ and $\beta$. A fan is a poset under inclusion, so we write $\alpha\prec\beta$ for $\alpha\subset\beta$ to emphasize this. If $\alpha\precneq\beta$ and $\alpha$ is maximal among cones contained in $\beta$, we write $\alpha\precdot\beta$. The set of codimension-$k$ cones in $\Delta$ we denote by $\Delta^{(k)}$. We use $\rho$ exclusively for rays (1-dimensional cones) and $\sigma$ exclusively for top-dimensional cones. 

The $T$-invariant affine open chart of $X$ corresponding to $\alpha$ is $U_\alpha:=\Spec\CC[\alpha^\vee\cap M]$, where $\alpha^\vee$ is the set of linear forms non-negative on $\alpha$.

There is an order reversing bijection between the $T$-orbit closures of $X$ and cones in $\Delta$. We write $O_\alpha$ for the $T$-orbit corresponding to $\alpha$, and denote its closure by $V(\alpha)$. In fact, $V(\alpha)$ is also a toric variety, with torus $T_\alpha$ the quotient of $T$ by the stabilizer of an element of $O_\alpha$. The character lattice for $T_\alpha$ is $M_\alpha=\alpha^\perp\subset M$. We denote the $\ZZ$-span of lattice points in $\alpha$ by $N^\alpha$, and the quotient $N_\alpha = N/N^\alpha\cong M_\alpha^\vee$ is the lattice of one parameter subgroups of $T_\alpha$. We use $\langle , \rangle$ to denote the pairing beween $M_\alpha$ and $N_\alpha$. When $\alpha\prec\beta$, the image of $\beta$ in $N_\alpha$ is a cone which we denote $\overline{\beta}$. If $\alpha\precdot\beta$, $\overline{\beta}$ is a ray, whose primitive generator we denote by $v_{\beta,\alpha}$, or $v_\beta$ if $\alpha=\{0\}$. 

Recall that for a simplicial cone $\beta$ with extremal rays generated by $v_1,\ldots,v_k$, its multiplicity is $\mult(\beta):=[N^\beta:\ZZ v_1+\ldots+\ZZ v_k]$. More generally for $\alpha\prec\beta$, we denote the multiplicity of the image of $\beta$ in $N_\alpha$ by $\mult_\alpha(\beta)$, so $\mult(\beta)=\mult_{\{0\}}(\beta)$ --geometrically, $\mult_\alpha(\beta)$ is the Hilbert-Samuel multiplicity of $U_\beta\cap V(\alpha)$ along $V(\beta)$. In the appendix, we show how to write a relative multiplicity $\mult_\alpha(\beta)$ in terms of usual multiplicities.

If $A$ is any abelian group, we denote the $\QQ$-vector space $A\otimes\QQ$ by $A_\QQ$.

\section{Background and generalities on Grothendieck groups and operational K-theory}

The Grothendieck group of coherent sheaves on an algebraic variety $X$ is denoted by $K_\circ(X)$. It is generated by isomorphism classes $[\mathscr{F}]$ for $\mathscr{F}$ a coherent sheaf on $X$, between which one imposes the relations $[\mathscr{F}]=[\mathscr{G}]+[\mathscr{H}]$ for each exact sequence $0\rightarrow\mathscr{G}\rightarrow\mathscr{F}\rightarrow\mathscr{H}\rightarrow 0$. For $f:Y\rightarrow X$ proper and $\mathscr{G}$ coherent on $Y$, there is a pushforward $f_*[\mathscr{G}]=\sum_i (-1)^i [R^i f_*\mathscr{G}]$. For $f:Y\rightarrow X$ flat and $\mathscr{F}$ coherent on $X$ there is a pullback $f^*[\mathscr{F}]=[f^*\mathscr{F}]$.

In $K$-theory the K\" unneth isomorphism does not hold in general. However, Anderson and Payne show in \cite[Proposition 6.4]{AP} that for $X$ toric it is true. In fact, it holds for \textit{linear varieties}, a larger class introduced by Totaro which includes toric varieties, see \cite{To} for the definition.

\begin{prop}
Let $X$ be a linear variety, and $Y$ an arbitrary variety. The natural map $K_\circ(X)\otimes_\ZZ K_\circ(Y)\rightarrow K_\circ(X\times Y)$ is an isomorphism.
\end{prop}

In what follows, we only consider tensor products of $K$-groups or their duals over $\ZZ$, so we omit the subscript. This theorem allows us to define a bilinear product on $K_\circ(X)^\vee:=\Hom_\ZZ(K_\circ(X),\ZZ)$ for $X$ any linear variety. Namely, let $\delta_*: K_\circ(X)\rightarrow K_\circ(X\times X)\cong K_\circ(X)\otimes K_\circ(X)$ be pushforward along the diagonal map. Then there is an induced map $\delta^*:K_\circ(X)^\vee\otimes K_\circ(X)^\vee\rightarrow K_\circ(X)^\vee$.

\begin{Defi}\label{def:proddef}
Let $c,d\in K_\circ(X)^\vee$. Set $c\cdot d := \delta^*(c\otimes d)$.
\end{Defi}

In simple terms, if $\delta_*(z)=\sum_i a_i\otimes b_i$, then $(c\cdot d)(z) = \sum_i c(a_i)d(b_i)$. Then it follows:

\begin{prop}\label{prop:isring}
For $X$ an arbitrary linear variety, $K_\circ(X)^\vee$ is a commutative ring.
\end{prop}
\noindent {\sl Proof:} 
The bilinear product defined is symmetric, since $\delta$ is equal to its composition with the involution on $X\times X$. The product is also associative, since $(\delta\times id)\circ\delta=(id\times \delta)\circ\delta$, so dually $\delta^*\circ(\delta^*\times id^*)=\delta^*\circ (id^*\times \delta^*)$. 
\qed
\\

Now, we verify that when $X$ is complete, the product we have defined on $K_\circ(X)^\vee$ is compatible with the one on operational $K$-theory. 

An element $c$ in $\opk^\circ(X)$ is a collection $(c_f)_{f}$ of endomorphisms of $K_\circ(Y)$ for each $f:Y\rightarrow X$. The collection $(c_f)_{f}$ must be compatible, in the sense that the maps must commute with proper pushforwards, flat pullbacks, and Gysin homomorphisms. Addition is defined by $(c_f)_f+(d_f)_f=(c_f+d_f)_f$, and the product, denoted $c\cup d$, is similarly given by $(c_f)_f\cup(d_f)_f=(c_f\circ d_f)_f$. For further details, we refer the reader to \cite[Section 4]{AP}. Amazingly, the product is commutative if $X$ admits a resolution of singularities (via the Kimura sequence \cite[Proposition 5.4]{AP}, though when $X$ is a complete linear variety it can also be deduced from the next proposition). We use the following shorthand: when there is a map $f:Y\rightarrow X$, and elements $z\in K_\circ(Y)$ and $c\in\opk^\circ(X)$, we drop the subscript $f$ and write $c\cap z$ for $c_f(z)$ when it does not cause confusion. When $X$ is complete, one may define a map from $\opk^\circ(X)$ to $K_\circ(X)^\vee$ sending $(c_f)_f$ to $\chi(c_{Id}(-))$. This map is an isomorphism if $X$ is also linear, a special case of \cite[Theorem 6.1]{AP}. 

Now, the main point of this section is the following, which is a K-theoretic analogue of \cite[Theorem 4]{FMSS}:

\begin{prop}
Let $X$ be a complete linear variety and let $\delta:X\rightarrow X\times X$ be the diagonal map. Given an expression $\delta_*(z)=\sum_i m_i a_i\otimes b_i$ with $m_i\in\QQ$, the product of classes $c$ and $d$ in $\opk^\circ(X)$ evaluated on $z$ satisfies:

\[
	\chi((c\cup d)\cap z)=\sum_i m_i \chi(c(a_i))\chi(d(b_i)).
\]

\end{prop}
\noindent {\sl Proof:}
For any morphism $f:Y\rightarrow X$, $c\in\opk^\circ(X)$, and $z\in K_\circ(Y)$, the identity 
\[
f^*c\cap z = \sum\chi(c\cap u_i)v_i, \tag{*}
\label{pullbackcap}
\]
holds, where $\gamma_f$ is the graph of $f$ and $(\gamma_f)_*(z) = \sum u_i\otimes v_i\in K_\circ(X)\otimes K_\circ(Y)$. To prove this, let $\pi_1$ and $\pi_2$ be the projections from $X\times Y$ to $X$ and $Y$. Then, one has the identities $\pi_2\circ\gamma_f=id_Y$, $\pi_1\circ\gamma_f=f$. Also, operational classes commute with proper pushfoward, so we have 
\[
f^*c\cap z = (id_Y)_*(f^*c\cap z) = (\pi_2)_*(\gamma_f)_*(\gamma_f^*\pi_1^*c\cap z) = (\pi_2)_*(\pi_1^*c\cap(\gamma_f)_*(z)).
\]
Substituting in our expression for $(\gamma_f)_*(z)$ and using the fact that flat pull-back and operational classes commute, we have
\[
(\pi_2)_*(\pi_1^*c\cap(\gamma_f)_*(z)) = (\pi_2)_*(\pi_1^*c\cap \sum u_i\otimes v_i) = (\pi_2)_*(\sum (c\cap u_i)\otimes v_i),
\]
and finally 
\[
(\pi_2)_*(\sum (c\cap u_i)\otimes v_i)=\sum(\pi_2)_*(c\cap u_i)\otimes v_i=\sum\chi(c\cap u_i)v_i,
\] 
since higher direct images commute with flat pull-back.

Then, in the context of the proposition we apply \eqref{pullbackcap} in the case that $f=\delta$. Since $d=\delta^*(id\otimes d)$, we obtain
\[
(c\cup d)\cap z= c\cap (d\cap \sum_i m_i a_i\otimes b_i) = c \cap (\sum_i m_i \chi(d(b_i))a_i) = \sum_i m_i \chi(d(b_i))c(a_i),
\]
to which we apply $\chi$ to obtain the proposition.
\qed

\begin{rem}\label{rem:robmring}
The results of this section \textcolor{black}{apply} mutatis mutandis when one replaces 
\begin{itemize}
	\item $K_\circ(X)^\vee$ with $E_*(X)^\vee:=\Hom_{E_*(pt)}(E_*(X),E_*(pt))$ for $E_*(-)$ a ``refined oriented Borel-Moore functor" with a K\" unneth isomorphism for some class of 	varieties, and
	\item $\opk^\circ(-)$ with the operational theory corresponding to $E_*(-)$.
\end{itemize}
When $X$ is linear, by results of \cite{FM,FMSS,To} this includes the pair $(A_*(X))^\vee$ and $A^*(X)$ where $A_*(X)$ is the Chow group. If we impose further that $X$ is toric this includes algebraic cobordism and its operational theory, for which the K\" unneth isomorphism was shown in \cite{A}, and the operational theory was developed in \cite{GK}.
\end{rem}

\section{K-theoretic balancing}

The definition of Grothendieck weights is quite simple:

\begin{Defi}
Let $\GW(\Delta)$ be the image of the map from $K_\circ(X)^\vee$ to $\ZZ^\Delta$ which sends $c$ to the function ($\alpha \rightarrow c([\mathcal{O}_{V(\alpha)}])$).
\end{Defi}

The map from $K_\circ(X)^\vee$ to $\ZZ^\Delta$ is injective since the classes $[\mathcal{O}_{V(\alpha)}]$ generate $K_\circ(X)$, so $\GW(\Delta)$ and $K_\circ(X)^\vee$ are isomorphic as abelian groups. Since we showed in Section 2 that if $X$ is linear then $K_\circ(X)^\vee$ is a ring, one can certainly pull back the ring structure to $\GW(\Delta)$. We would like to consider now how to characterize $\GW(\Delta)$ inside $\ZZ^\Delta$, and whether it is possible to give a displacement rule for products in $\GW(\Delta)$ like the one for Minkowski weights. 

The first point to note is that characterizing $\GW(\Delta)$ inside all $\ZZ$-valued functions on $\Delta$ is equivalent to characterizing the relations between $[\mathcal{O}_{V(\alpha)}]$. Let us introduce some notation for the next proposition: If $P$ is a lattice polytope in $M$ with normal fan $\Delta$ in $N$, we write $P\perp\Delta$. Recall that there is a containment-reversing correspondence between the cones of $\Delta$ and the faces of $P$, which are themselves lattice polytopes. For $\alpha\in\Delta$, let the corresponding face of $P$ be denoted $F_\alpha$. Also, we recall that the Ehrhart polynomial $L(P,t)$ of $P$ is the polynomial determined by $L(P,t_0)=|t_0P\cap M|$.

\begin{prop}\label{prop:ehrhartequiv}
Let $\Delta$ be a projective simplicial fan and $X$ the corresponding toric variety. Let $\phi_P:\QQ^\Delta\rightarrow\QQ[t]$ be the map sending the tuple $(a_\alpha)$ to $\sum_{\alpha\in\Delta}a_\alpha L(F_\alpha,t)$. Let $\psi:\QQ^\Delta\rightarrow K_\circ(X)_\QQ$ be the map sending $e_\alpha$ to $[\mathcal{O}_{V(\alpha)}]$. Then
\[
\ker(\psi)=\bigcap_{P\perp\Delta}\ker(\phi_P).
\]
\end{prop}
\noindent {\sl Proof:}
The ``$\subset$" direction is a translation of the well-known vanishing of higher cohomology for ample line bundles on toric varieties: let $\sum_{\alpha\in\Delta}a_\alpha[\mathcal{O}_{V(\alpha)}]=0$ in $K_\circ(X)_\QQ$, and for $P$ any lattice polytope with normal fan $\Delta$ let $L_{t_0P}, t_0\in\ZZ_{\geq 0}$ be the ample line bundle that corresponds to the $t_0$-th dilate. Then
\[
0=\chi([L_{t_0P}]\cdot(\sum_{\alpha\in\Delta}a_\alpha[\mathcal{O}_{V(\alpha)}])) = \sum_{\alpha\in\Delta}a_\alpha \chi([L_{t_0P}|_{V(\alpha)}])=\sum_{\alpha\in\Delta}a_\alpha L(F_\alpha,t_0).
\]
Thus the polynomial $\sum_{\alpha\in\Delta}a_\alpha L(P,t)$ has infinitely many roots, and so must be $0$. 

Now we verify the other direction: Suppose, for any polytope $P$, that \\ $\sum_{\alpha\in\Delta}a_\alpha L(F_\alpha,t)=0$. Translating to geometry, this states (when $t=1$) that
\[
0=\sum_{\alpha\in\Delta}a_\alpha\chi(L_P|_{V(\alpha)})=\chi\left(\left(\sum_{\alpha\in\Delta}a_\alpha[\mathcal{O}_{V(\alpha)}]\right)\cdot [L_P]\right).
\]
Then, the result follows from the next lemma. In the proof, we write  ``$\nicedeg$" for the projection map from $A^*(X)$ to $A^0(X)$.

\begin{lem}
Suppose $\Delta$ is a simplicial projective fan and $X$ the corresponding toric variety. If for $x\in K_\circ(X)$, we have $\chi(x\cdot [L_P])=0$ for all $P$ with normal fan $\Delta$, then $x=0$ in $K_\circ(X)_\QQ$.
\end{lem}

\noindent {\sl Proof:} 
To avoid clutter, we assume in this proof that all Chow and $K$-groups and their duals are tensored with $\QQ$. When $X$ is projective it is well-known that ample divisors generate the Neron-Severi group of divisors modulo numerical equivalence $NS(X)$ (which is the same as $A_{n-1}(X)$ if $X$ is toric). Additionally, on a complete toric variety $A^*(X)\cong A_*(X)^\vee$ via $c\rightarrow \nicedeg(c\cap-)$, by \cite[Theorem 3]{FMSS}. Thus, if for all $c\in A^*(X)$, $\nicedeg(c\cap y)=0$, then $y$ must be zero. Since $X$ is simplicial, $A^*(X)$ is generated as an algebra by Chern classes of $T$-equivariant divisors, which are in turn generated as a group by $T$-equivariant ample line bundles. Thus, $A^*(X)$ is generated over $\QQ$ by $1$ and monomials in $c_1(L_P)$ for $P$ a polytope. Equivalently, $A^*(X)$ is generated by $1$ and $1+c_1(L_P)=\chern(L_P)$ for $\chern:K^\circ(X)\rightarrow A^*(X)$ the Chern character map. Since $\chern$ is a ring homomorphism, $\chern(L_P)\chern(L_Q)=\chern(L_P\otimes L_Q)$. But $L_P\otimes L_Q=L_{P+Q}=L_{P'}$. Thus, $y$ is zero if $\nicedeg(y)=0$ and
\[
\nicedeg(\chern(L_{Po})\cap y)=0,
\]
for all $P$.
By Riemann-Roch for algebraic schemes as in \cite[Chapter 18]{F98},
\[
\nicedeg(\chern(L_P)\cap y)=\nicedeg(\tau_X(L_P\otimes\tau_X^{-1}(y)))=\chi(L_P\otimes\tau_X^{-1}(y)).
\]
Since $\tau_X$ is an isomorphism between $K_\circ(X)$ and $A_*(X)$, the lemma and proposition are proved.
\qed
\\

Returning to the question of an explicit balancing condition, Proposition \ref{prop:ehrhartequiv} indicates that for an arbitrary fan, directly finding a simple set of generators for $\ker(\psi)$ (equivalently a combinatorial characterization for Grothendieck weights) may be difficult, in comparison to the situation for $\MW^*(\Delta)$. However, using the Riemann-Roch transformation gives us an avenue to use balancing conditions for Minkowski weights to give ones for Grothendieck weights. Recall that this is a map $\tau_X:K_\circ(X)\rightarrow A_*(X)_\QQ$, which is an isomorphism after tensoring $K_\circ(X)$ with $\QQ$. The map commutes with proper pushforwards, and for a vector bundle $E$, there is the equality $\tau_X([E])=\chern(E)\cdot td(X)$. For details regarding the Todd class $td(X)$ for singular varieties and this version of the Riemann-Roch theorem, see \cite[Chapter 18]{F98}.

To make things concrete, we first suppose that $\Delta$ is simplicial. In this case Pommersheim and Thomas introduced in \cite{PT} certain rational numbers $t^\alpha_\rho$ for each cone $\alpha$ and ray $\rho$ such that $\rho\subset\alpha$, which depend on the choice of a generic complete flag $F_\bullet$ in $N_\QQ$. 

\begin{Defi}\label{def:tref}
Let $F_\bullet$ be a generic complete flag in $N_\QQ$, so $F_i$ is an $i$-dimensional subspace of $N_\QQ$. Given $\alpha\in\Delta$ of dimension $k$, let the primitive elements of the rays of $\alpha$ be $v_{\rho_1},\ldots,v_{\rho_k}$. Then by genericity, $F_{n-k+1}\cap\QQ\cdot\alpha$ is $1$-dimensional, and so determines a vector (unique up to scaling):
\[
0\neq\sum_{i=1}^k t_{\rho_i}^\alpha v_{\rho_i}\in F_{n-k+1}\cap \QQ\cdot\alpha,
\]
We only consider generic $F_\bullet$ such that all $t_\rho^\alpha$ are non-zero.
\end{Defi}

We require these $t^\alpha_\rho$ for the following result from \cite{PT}: an explicit formula for monomials in the $T$-invariant divisors $[V(\rho)]$ as a $\QQ$-linear combination of classes of subvarieties $[V(\alpha)]$. For a cone $\beta$, let $\Rho_\beta$ be the set of rays in $\beta$. Let $S$ be some set of rays in $\Delta$. For $\rho\in S$, let $a_\rho$ be some positive integers such that $\sum_{\rho\in S}a_\rho = l$. Then a restatement of their Theorem 3 is:
\begin{prop}\label{prop:absolutedivisorsproduct}
\[
\prod_{\rho\in S} [V(\rho)]^{a_\rho}=\sum_{\substack{\beta\in\Delta^{(n-l)}\\S\subset\Rho_\beta}}\frac{\prod_{\rho\in S}(t_\rho^\beta)^{a_\rho}}{\prod_{\rho\in\Rho_\beta}t_\rho^\beta}\frac{[V(\beta)]}{\mult(\beta)}.
\]
\end{prop}
For $\alpha\in\Delta$ and a generic flag $F_\bullet$ in $N$, the images $\overline{F_1}\subset\ldots\subset\overline{F_{n-\dim(\alpha)}}$ in $N_\alpha$ form a generic flag. Thus for $\beta$ a cone containing $\alpha$ and $\rho\in\Rho_\beta\smallsetminus\Rho_\alpha$ there are also numbers $t^{\overline{\beta}}_{\overline{\rho}}$. We relate these to $t_\rho^\beta$ in the next proposition.

\begin{prop}\label{prop:multratio}
Let $\alpha\prec\beta$ be simplicial cones in a fan $\Delta$. Then for $\rho\in\Rho_\beta\smallsetminus\Rho_\alpha$, $t^{\overline{\beta}}_{\overline{\rho}}=\frac{\mult(\alpha+\rho)}{\mult(\alpha)}t^\beta_\rho$.
\end{prop}
\noindent {\sl Proof:} 
The unique vector in $\overline{F_{n-k+1}}\cap\QQ\cdot\overline{\beta}$ is the image of the unique vector in $F_{n-k+1}\cap\QQ\cdot\beta$, which in explicit terms is
\[
\overline{\sum_{\rho\in\Rho_\beta} t^\beta_{\rho}v_{\rho}}=\sum_{\rho\in\Rho_\beta\smallsetminus\Rho_\alpha} t_{\rho}^\beta\overline{v_{\rho}}.
\]
But the image of a primitive generator of a ray $v_\rho$ is not necessarily primitive, i.e. $\overline{v_\rho}=b_\rho v_{\overline{\rho}}$ for $b_\rho$ a positive integer. In fact, $b_\rho=[\ZZ v_{\overline{\rho}}:\ZZ \overline{v_{\rho}}]$ of subgroups of $N_\alpha$. If $\pi_\alpha:N\rightarrow N_\alpha$ is the quotient map, then $\pi_\alpha^{-1}(\ZZ v_{\overline{\rho}}) = N^{\alpha+\rho}$, and $\pi_\alpha^{-1}(\ZZ \overline{v_{\rho}})=N^\alpha+\ZZ v_\rho$. Thus $[\ZZ v_{\overline{\rho}}:\ZZ \overline{v_{\rho}}]=[N^{\alpha+\rho}:N^\alpha+\ZZ v_\rho]$. Then we can decompose $\mult(\alpha+\rho)$ as a product:
\begin{align*}
\mult(\alpha+\rho) 	& =[N^{\alpha+\rho}:\ZZ v_1+\ldots+\ZZ v_k+\ZZ v_\rho]\\
				& = [N^{\alpha+\rho}:N^\alpha+\ZZ v_\rho][N^\alpha+\ZZ v_\rho:\ZZ v_1+\ldots+\ZZ v_k+\ZZ v_\rho].
\end{align*}
But $[N^\alpha+\ZZ v_\rho:\ZZ v_1+\ldots+\ZZ v_k+\ZZ v_\rho]= [N^\alpha:\ZZ v_1+\ldots \ZZ v_k]=\mult(\alpha)$. Thus $[N^{\alpha+\rho}:N^\alpha+\ZZ v_\rho]=\frac{\mult({\alpha+\rho})}{\mult(\alpha)}=b_\rho$. Thus, we have
\[
\sum_{\rho\in\Rho_\beta\smallsetminus\Rho_\alpha} t_{\rho}^\beta\overline{v_{\rho}} = \sum_{\rho\in\Rho_\beta\smallsetminus\Rho_\alpha} t_{\rho}^\beta\frac{\mult(\alpha+\rho)}{\mult(\alpha)}v_{\overline{\rho}}.
\]
Since the primitive generators of the rays in $\overline{\beta}$ are the $v_{\overline{\rho}}$, we are done by the definition of $t^{\overline{\beta}}_{\overline{\rho}}$.
\qed
\\

For each cone $\alpha\in\Delta$ Brion and Vergne defined the (finite) subgroup $G_\alpha\subset(\CC^*)^{\Rho_\alpha}$ to be the kernel of the map $(\CC^*)^{\Rho_\alpha}\rightarrow T$ given by 
\[
(c_\rho)_\rho \mapsto \prod_{\rho\in\Rho_\alpha}v_\rho(c_\rho).
\]
For nested cones $\alpha\prec\beta$, we use the notation $G^\alpha_\beta$ for the analogous subgroup defined from the data of $\overline{\beta}$ in the quotient fan. Explicitly, this is the kernel of the map $(\CC^*)^{\Rho_\beta\smallsetminus\Rho_\alpha}\rightarrow T_\alpha$ given by
\[
(c_\rho)_\rho \mapsto \prod_{\rho\in\Rho_\beta\smallsetminus\Rho_\alpha}v_{\overline{\rho}}(c_\rho).
\]
Let the number of rays in $\Delta_\alpha$ be $k$. We define $G_{\Delta_\alpha}$ to be the union inside $(\CC^*)^{k}$ of $G^\alpha_\beta$ over all $\beta$ containing $\alpha$. For $\rho$ a ray in $\beta\smallsetminus\alpha$, we denote the character $G^\alpha_\beta\rightarrow\CC^*$ given by projection by $a^\alpha_{\rho}$.

Using these numbers, we have the following proposition. We will refer to the degree $0$ term of a formal Laurent series $f(t_1,\ldots,t_k)$ by $f(t_1,\ldots,t_k)_{[0]}$, i.e. if 
\[f=\left(\frac{1}{1-e^{-t}}\right)\left(\frac{1}{1-e^{-s}}\right)=\left(\frac{1}{t}+\frac{1}{2}+\frac{t}{12}+\Oc(t^2)\right)\left(\frac{1}{s}+\frac{1}{2}+\frac{s}{12}+\Oc(s^2)\right),\]
then 
\[
f_{[0]}=\frac{1}{4}+\frac{1}{12}\left(\frac{t}{s}+\frac{s}{t}\right).
\]

\begin{prop}\label{prop:rrmtx}
For $X$ a complete  simplicial toric variety, the Riemann-Roch transformation has the form 
\[
\tau_X([\mathcal{O}_{V(\alpha)}]) = \sum_{\alpha\prec\beta}\sum_{g\in G^\alpha_\beta}\left(\prod_{\rho\in\Rho_\beta\smallsetminus\Rho_\alpha}\frac{\mult(\alpha+\rho)/\mult(\alpha)}{1-a^\alpha_{\rho}(g)e^{-mult(\alpha+\rho)t^\beta_\rho}}\right)_{[0]}\frac{[V(\beta)]}{\mult(\beta)}.
\]
\end{prop}
\noindent {\sl Proof:} 
In section 4.2 of \cite{BV}, the authors provide a formula for the Todd class of a complete simplicial toric variety. Applied to $V(\alpha)$, this gives:
\[
\tau_{V(\alpha)}([\mathcal{O}_{V(\alpha)}]) = \sum_{g\in G_{\Delta_\alpha}} \prod_{\rho\in\Delta_\alpha^{(n-\dim(\alpha)-1)}}\frac{[V(\rho)]}{1-a^\alpha_\rho(g)e^{-[V(\rho)]}}.
\]
Since the Todd class commutes with proper pushfoward, we have 
\[
\tau_X([\mathcal{O}_{V(\alpha)}]) = \sum_{g\in G_{\Delta_\alpha}} \prod_{\rho\in\Delta_\alpha^{(n-\dim(\alpha)-1)}}i_*\frac{[V(\rho)]}{1-a^\alpha_\rho(g)e^{-[V(\rho)]}}.
\]
Applying Proposition \ref{prop:absolutedivisorsproduct}, we have
\begin{align*}
i_*\prod_{\rho\in\Delta^{(n-\dim(\alpha)-1)}_\alpha}\frac{[V(\rho)]}{1-a^\alpha_{\rho}(g)e^{-[V(\rho)]}} &=\sum_{\alpha\prec\beta}\left(\prod_{\rho\in\Rho_\beta\smallsetminus\Rho_\alpha} \frac{1}{1-a^\alpha_{\rho}(g)e^{-t_{\rho}^{\overline{\beta}}}} \right)_{[0]}\frac{i_*[V(\overline{\beta})]}{\mult_\alpha(\beta)},\\
&=\sum_{\alpha\prec\beta}\left(\prod_{\rho\in\Rho_\beta\smallsetminus\Rho_\alpha} \frac{1}{1-a^\alpha_{\rho}(g)e^{-t_{\overline{\rho}}^{\overline{\beta}}}} \right)_{[0]}\frac{[V(\beta)]}{\mult_\alpha(\beta)}.
\end{align*}
Using the formula for $t_{\overline{\rho}}^{\overline{\beta}}$ from the previous proposition, and the formula for $\mult_\alpha(\beta)$ from the appendix, we obtain that the above is equal to 
\[
\sum_{\alpha\prec\beta}\left(\prod_{\rho\in\Rho_\beta\smallsetminus\Rho_\alpha} \frac{1}{1-a^\alpha_{\rho}(g)e^{-\frac{\mult(\alpha+\rho)}{\mult(\alpha)}t_{\rho}^{\beta}}} \right)_{[0]}\frac{[V(\beta)]}{\mult(\beta)\prod_{\Rho_\beta\smallsetminus\Rho_\alpha}\frac{\mult(\alpha)}{\mult(\alpha+\rho)}}.
\]
Due to the ``degree 0'' imposition, the $\mult(\alpha)$ factor in the exponent $e^{-\frac{\mult(\alpha+\rho)}{\mult(\alpha)}t^\beta_\rho}$ can be cancelled. Summing over $g\in G^\alpha_\beta$ gives the proposition.
\qed

\begin{ex}\label{ex:rrmtx}
We use this proposition to calculate the $t^\alpha_\rho$ and Riemann-Roch matrix for a weighted projective space $X:=\PP(1,1,2,3)$. Recall that the fan of $X$ has rays $\rho_1=(1,0,0),\rho_2=(0,1,0),\rho_3=(0,0,1),$ and $\rho_4=(-1,-2,-3)$. The maximal cones are those generated by $3$-element subsets of $\{\rho_1,\rho_2,\rho_3,\rho_4\}$. If our flag in $\QQ^3$ is is given by
\[
\{0\}\subsetneq \nicespan\{(a,b,c)\}\subsetneq \nicespan\{(a,b,c),(d,e,f)\}\subsetneq\QQ^3,
\]
then we have the table of expressions for $t^\alpha_\rho$ in Table \ref{tab:table1}.

\begin{table}[h!]
  \begin{center}
    \caption{Example \ref{ex:rrmtx}}
    \label{tab:table1}
    \begin{tabular}{l|c|r} 
    \toprule
      Cone $(\alpha)$ & Ray $(\rho)$ & $t^\sigma_\rho$ \\
      \midrule
$\sigma_{123}$ & $\rho_1$ & $a$ \\
& $\rho_2$ & $b$ \\
& $\rho_3$ & $c$ \\
\hline
$\sigma_{124}$ & $\rho_1$ & $a-c /3$ \\
& $\rho_2$ & $b-2c/3$ \\
& $\rho_4$ & $-c/3$ \\
\hline
$\sigma_{134}$ & $\rho_1$ & $a-b / 2$ \\
& $\rho_3$ & $c-3b/ 2$ \\
& $\rho_4$ & $-b/ 2$ \\
\hline
$\sigma_{234}$ & $\rho_2$ & $b-2a$ \\
& $\rho_3$ & $c-3a$ \\
& $\rho_4$ & $-a$ \\
\hline
$\alpha_{12}$ & $\rho_1$ & $af-cd$ \\
& $\rho_2$ & $bf-ce$ \\
\hline
$\alpha_{13}$ & $\rho_1$ & $ae-bd$ \\
& $\rho_3$ & $ce-bf$ \\ 
\hline
$\alpha_{14}$ & $\rho_1$ & $3(ae-bd)+ 2(cd-af)+(bf-ce)$ \\
& $\rho_4$ & $bf-ce$ \\
\hline
$\alpha_{23}$ & $\rho_2$ & $bd-ae$ \\
& $\rho_3$ & $af-cd$ \\
\hline
$\alpha_{24}$ & $\rho_2$ & $-(3(ae-bd)+ 2(cd-af)+(bf-ce))$ \\
& $\rho_4$ & $af-cd$ \\
\hline
$\alpha_{34}$ & $\rho_3$ & $3(ae-bd)+ 2(cd-af)+(bf-ce)$ \\
& $\rho_4$ & $ae-bd$ \\
\hline
$\rho$ & $\rho$ & $1$ \\
\bottomrule
    \end{tabular}
  \end{center}
\end{table}
\end{ex}

After calculating these $t^\alpha_\rho$, one can write the Todd class of each subvariety in a uniform way with rational functions the $t_\rho^\sigma$ as coefficients. The column vector corresponding to the image of $[\mathcal{O}_X]$ is \textcolor{black}{on the last page}. For the flag specified by $(a,b,c)=(2,3,5)$, $(d,e,f)=(3,5,7)$, one obtains the Riemann-Roch matrix in Table \ref{tab:table2}.

\begin{table}[h!]
\footnotesize
  \begin{center}
    \caption{The Riemann-Roch matrix of Example \ref{ex:rrmtx}}
    \label{tab:table2}
    \tabcolsep=0.11cm
    \begin{tabular}{c|c|c|c|c|c|c|c|c|c|c|c|c|c|c|c} 
    \toprule
\tiny{${}_{[V(-)]}\setminus^{[\mathcal{O}_{V(-)}]}$} & $X$ & $\rho_1$ & $\rho_2$ & $\rho_3$ & $\rho_4$ & $\alpha_{12}$ & $\alpha_{13}$ & $\alpha_{14}$ & $\alpha_{23}$ & $\alpha_{24}$ & $\alpha_{34}$ & $\sigma_{123}$ & $\sigma_{124}$ & $\sigma_{134}$ & $\sigma_{234}$\\
      \midrule
$X$ &1 & 0 & 0& 0& 0& 0& 0& 0& 0& 0& 0& 0& 0& 0& 0\\
$\rho_1$&1/2 & 1 & 0& 0& 0& 0& 0& 0& 0& 0& 0& 0& 0& 0& 0\\
$\rho_2$&1/2 & 0& 1& 0& 0& 0& 0& 0& 0& 0& 0& 0& 0& 0& 0\\
$\rho_3$&1/2 & 0& 0& 1& 0& 0& 0& 0& 0& 0& 0& 0& 0& 0& 0\\
$\rho_4$&1/2 & 0& 0& 0& 1& 0& 0& 0& 0& 0& 0& 0& 0& 0& 0\\
$\alpha_{12}$&29/48 & 1/2& 1/2& 0& 0& 1& 0& 0& 0& 0& 0& 0& 0& 0& 0\\
$\alpha_{13}$&29/48 & 1/2& 0& 1/2& 0& 0& 1& 0& 0& 0& 0& 0& 0& 0& 0\\
$\alpha_{14}$&-5/48 & 1/2& 0& 0& 1/2& 0& 0& 1& 0& 0& 0& 0& 0& 0& 0\\
$\alpha_{23}$&5/12 & 0& 1/2& 1/2& 0& 0& 0& 0& 1& 0& 0& 0& 0& 0& 0\\
$\alpha_{24}$&5/12 & 0& 1/2& 0& 1/2& 0& 0& 0& 0& 1& 0& 0& 0& 0& 0\\
$\alpha_{34}$&5/12 & 0& 0& 1/2& 1/2& 0& 0& 0& 0& 0& 1& 0& 0& 0& 0\\

$\sigma_{123}$&31/72  & 79/180 & 59/120 & 3/2 & 0& 1/2& 1/2& 0& 1/2& 0& 0& 1& 0& 0& 0\\

$\sigma_{124}$&1/8 & 41/60 & -11/60 & 0& 1/12 & 1/2& 0& 1/2& 0& 1/2& 0& 0& 1& 0& 0\\
    
$\sigma_{134}$&1/36 & 1/9& 0 & 1/9 & 1/3& 1/2& 1/2& 0& 0& 1/2& 0 & 0 & 0 & 1 & 0\\
    
$\sigma_{234}$&5/12 & 0 & 11/24 & 11/24 &  5/12 & 0& 0& 0& 1/2& 1/2& 1/2& 0& 0& 0& 1\\
\bottomrule
    \end{tabular}
  \end{center}
\end{table}

\begin{Defi}\label{def:invcoeffs}
Let $\mu_\alpha(\beta)$ be $\sum_{g\in G^\alpha_\beta}\left(\prod_{\rho\in\Rho_\beta\smallsetminus\Rho_\alpha}\frac{\mult(\alpha+\rho)/\mult(\alpha)}{1-a^\alpha_{\rho}(g)e^{-mult(\alpha+\rho)t^\beta_\rho}}\right)_{[0]}$. Then, the matrix $(\mu_\alpha(\beta))$ defines an isomorphism from $\QQ^\Delta$ to itself lifting the Riemann-Roch isomorphism from $K_\circ(X)_\QQ$ to $A_*(X)_\QQ$. Let $\nu_\alpha(\beta)$ refer to the $(\alpha,\beta)$-th entry of the inverse of $(\mu_\alpha(\beta))$, so $\tau_X^{-1}([V(\alpha)]) = \sum_{\alpha\prec\beta}\nu_\alpha(\beta) [\mathcal{O}_{V(\beta)}]$. Then for example $\nu_\alpha(\alpha)=1$, and for $\Delta$ smooth and the $t_\rho^\alpha$ as defined in \ref{def:tref} we can write:
\[
\nu_{\alpha}(\beta) = \sum_{\substack{\alpha=\alpha_0\precneq\alpha_1\\ \precneq\ldots\precneq\alpha_k=\beta}}(-1)^k\left(\prod_{l=1}^k\left(\prod_{\rho\in\alpha_l\smallsetminus\alpha_{l-1}}\frac{1}{1-e^{-t^{\alpha_l}_\rho}}\right)_{[0]}\right).
\]
\end{Defi} 

Finally, we can prove:

\begin{thm}\label{thm:balancing}
For $\Delta$ a complete simplicial fan, a function $g:\Delta\rightarrow\ZZ$ is a Grothendieck weight if and only if it satisfies
\[
	\sum_{\alpha\precdot\beta}\langle u, v_{\beta,\alpha}\rangle \sum_{\beta\prec\gamma}\nu_{\beta}(\gamma) g(\gamma) = 0,
\]
for all $\alpha\in\Delta, u\in M(\alpha)$.
\end{thm}

\noindent {\sl Proof:} 
The result \cite[Proposition 2.1]{FS} says that the expressions $\sum_{\alpha\precdot\beta}\langle u, v_{\beta,\alpha}\rangle [V(\beta)]$ generate the kernel of the map $\ZZ^{\Delta}\rightarrow A_*(X)$ sending $e_\alpha$ to $[V(\alpha)]$. One thus has an isomorphism of exact sequences:
\[
\begin{tikzcd}
0 \arrow[r]& \Rel_{K_\QQ} \arrow[r]\arrow[d]& \QQ^\Delta \arrow[r]\arrow[d]& K_\circ(X)_\QQ \arrow[r]\arrow[d]& 0 \\
0 \arrow[r]& \Rel_{A_\QQ} \arrow[r] & \QQ^\Delta \arrow[r] & A_*(X)_\QQ \arrow[r] & 0.
\end{tikzcd}
\]
where the map on the right is $\tau_X$, and the middle and left maps are given by sending $e_\alpha$ to $\sum_{\alpha\prec\beta}\mu_\alpha(\beta) e_\beta$. The inverse image of $\sum_{\alpha\precdot\beta}\langle u, v_{\beta,\alpha}\rangle [V(\beta)]\in \Rel_A$ is $\sum_{\alpha\precdot\beta}\langle u, v_{\beta,\alpha}\rangle \sum_{\beta\prec\gamma}\nu_\beta(\gamma)[\mathcal{O}_{V(\gamma)}]$, so such relations generate $\Rel_{K_\QQ}$, and appropriate multiples of these relations generate a finite index subgroup of $\Rel_K$. Dually,  $K_\circ(X)^\vee$ must then consist of linear forms sending such expressions to $0$, which implies that the relations in the theorem statement characterize Grothendieck weights.
\qed
\\

The following lemma explains how one can approach non-simplicial fans. It follows directly from Lemma \ref{lem:stability}, which we have kept in the last section on tropical geometry due to its relevance there.
\begin{lem}\label{lem:refinementlemma}
Let $\Delta'$ be an arbitrary fan and $\Delta$ a smooth refinement. Then $g:\Delta'\rightarrow\ZZ$ is a Grothendieck weight if and only if the function on $\Delta$ determined by $\alpha\rightarrow g(\alpha')$ for $\alpha'$ the smallest cone in $\Delta'$ containing $\alpha$ is a Grothendieck weight on $\Delta$.
\end{lem}

\begin{rem}
We could have used a different set of generators for the Grothendieck group in our definition of Grothendieck weights, e.g. ideal sheaves or canonical sheaves of invariant subvarieties. However, the problem of combinatorially describing the relations between these classes seems equally difficult.
\end{rem}

\begin{ex}
The characterization of Grothendieck weights given in Theorem \ref{thm:balancing} boils down to the following straightforward criteria for arbitrary complete fans up to dimension 3:
\begin{enumerate}
	\item Dimension 1: A $\ZZ$-valued function on the fan of $\PP^1$ is a Grothendieck weight if and only if it has the same value on 		both maximal cones.
	\item Dimension 2: A $\ZZ$-valued function $g$ on $\Delta$ the fan of a complete toric surface is a Grothendieck weight if and only if it is 		constant on maximal cones, and 
		\[ \sum_{\rho\in\Delta^{(1)}}\left(g(\rho)-g(\sigma)\right) v_\rho = 0,\]
		for $\sigma$ any maximal cone.
	\item Dimension 3: A $\ZZ$-valued function $g$ on $\Delta$ the fan of a complete toric threefold is a Grothendieck weight if and only if it is 	constant on maximal cones and, still writing $\sigma$ for any maximal cone:
		\begin{enumerate}
			\item $\sum_{\rho\precdot\beta}\left(g(\beta)-g(\sigma)\right) v_{\beta,\rho} = 0$
				for any fixed ray $\rho$, and
			\item $\sum_{\rho\in\Delta^{(2)}}\left(g(\rho)-\sum_{\rho\precdot\alpha}\frac{g(\alpha)}{2}\right)v_\rho=g(\sigma)						\left(\sum_{\rho\in\Delta^{(2)}}(1-\sum_{\rho\precdot\alpha}\frac{1}{2})v_\rho\right)$.
		\end{enumerate}
\end{enumerate}
\end{ex}

\section{Products of Grothendieck weights}

Let $\Delta$ be a complete simplicial fan, and suppose that $g_1$ and $g_2$ in $\GW(\Delta)$ are given. Their product, as induced by the product on $K_\circ(X)^\vee$ given in Definition \ref{def:proddef}, may be calculated explicitly via the formula in the next theorem. As in the case of Minkowski weights, to undertake the calculation one must choose an auxiliary ``displacement" vector in $N$ which we call $v$. Then for three cones $\alpha$, $\beta$, and $\gamma$, we define $m^\alpha_{\beta,\gamma}$ in the same manner as \cite{FS}, by

\begin{equation*}
    m^\alpha_{\beta,\gamma} = \begin{cases}
	       0					& \text{if }\textrm{$\beta\cap(\gamma+v)=\emptyset$,}\\
[N:\ZZ\cdot\beta+\ZZ\cdot\gamma] 	& \text{otherwise.}
           \end{cases}
\end{equation*}
Though it is suppressed in the notation, we emphasize that $m_{\beta,\gamma}^\alpha$ depends on the choice of $v$. Recall that for $f$ a formal Laurent series, $f_{[0]}$ denotes the degree $0$ term of $f$.

\begin{thm}\label{thm:prod}
The product weight $g_3 = g_1\cup g_2$ is given by
\[
g_3(\alpha) = \sum_{\alpha\prec\beta}\mu_\alpha(\beta)\sum_{\substack{\beta\prec\gamma,\epsilon \\ \codim(\gamma)+\codim(\epsilon)\\=\codim(\beta)}}m^\beta_{\gamma,\epsilon}\sum_{\substack{\gamma\prec\zeta\\ \epsilon\prec\eta}}\nu_{\gamma}(\zeta)\nu_{\epsilon}(\eta)g_1(\zeta)g_2(\eta),
\]
where $\mu_\alpha(\beta), \nu_{\alpha}(\beta)$ was defined in \ref{def:invcoeffs}.
\end{thm}
\noindent {\sl Proof:} 
Recall that $\delta:X\rightarrow X\times X$ is the diagonal map. Then, we have
\begin{align*}
\delta_*([\Oc_{V(\alpha)}]) & = \delta_*(\tau_X^{-1}(\tau_X([\Oc_{V(\alpha)}]))),\\
& = \tau_X^{-1}\left(\sum_{\alpha\prec\beta}\mu_\alpha(\beta)(\delta_*([V(\beta)]))\right).
\end{align*}
By \cite[Theorem 4.2]{FS}, we may use the $m^\beta_{\gamma,\epsilon}$ determined by our generic displacement vector to decompose each $\delta_*([V(\beta)])$, thus obtaining
\begin{align*}
& \tau_X^{-1}\left(\sum_{\alpha\prec\beta}\mu_\alpha(\beta)\delta_*([V(\beta)])\right),\\
= & \tau_X^{-1}\left(\sum_{\alpha\prec\beta}\mu_\alpha(\beta)\sum_{\substack{\beta\prec\gamma,\epsilon \\ \codim(\gamma)+\codim(\epsilon)\\=\codim(\beta)}} m^\beta_{\gamma,\epsilon}[V(\gamma)\times V(\epsilon)]\right),\\
= & \sum_{\alpha\prec\beta}\mu_\alpha(\beta)\sum_{\substack{\beta\prec\gamma,\epsilon \\ \codim(\gamma)+\codim(\epsilon)\\=\codim(\beta)}} m^\beta_{\gamma,\epsilon}\tau_X^{-1}([V(\gamma)\times V(\epsilon)]).
\end{align*}
But then,
\begin{align*}
\tau^{-1}([V(\gamma)\times V(\epsilon)]) = 	& \tau^{-1}([V(\gamma)])\tau^{-1}([V(\epsilon)]), \\
							 = 	& \sum_{\substack{\gamma\prec\zeta \\ \epsilon\prec\eta} } \nu_\gamma(\zeta)\nu_\epsilon(\eta)[\mathcal{O}_{V(\zeta) \times V(\eta)}].
\end{align*}
So, by definition (\ref{def:proddef}) of the product the theorem follows.
\qed
\\

We can use this theorem to show the following proposition and its corollary, which give some structure to $\GW(\Delta)$.

\begin{prop}\label{prop:vanishingideal}
Let $\Sigma\subset\Delta$ be fans. The set of Grothendieck weights on $\Delta$ that vanish on the complement of $\Sigma$ forms an ideal in $\GW(\Delta)$.
\end{prop}
\noindent {\sl Proof:} 
Suppose we have two weights $g_1,g_2$, such that $g_1$ vanishes on the cones of $\Delta$. Then in general there are some coefficients $C$ such that  
\[
(g_1\cup g_2) (\alpha) = \sum_{\alpha\prec\beta}\sum_{\substack{\beta\prec\gamma,\epsilon \\ \codim(\gamma)+\codim(\epsilon) \\ = \codim(\beta)}}\sum_{\substack{\gamma\prec\zeta \\ \epsilon\prec\eta}}C_{\substack{\alpha,\beta, \gamma, \epsilon, \\ \zeta,\eta}} g_1(\zeta)g_2(\eta).
\]
If $\alpha$ is not in $\Sigma$, then since $\alpha\prec\zeta$ and $\Sigma$ is a fan, certainly $\zeta$ is not in $\Sigma$. Thus $g_1(\zeta)$ must be $0$ for each term in the sum.
\qed

\begin{cor}\label{cor:filtration}
The ring $\GW(\Delta)$ is filtered by ideals $I_k$ consisting of functions that vanish on cones of codimension less than $k$.
\end{cor}

\begin{ex}
We calculate the product of the following Grothendieck weights:
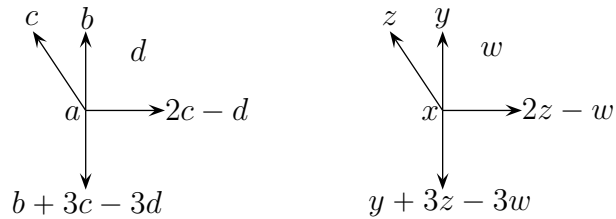
\[
\begin{pspicture}(-35,-50)(100,50)

\psline{->}(0,0)(0,30)
\psline{->}(0,0)(0,-30)
\psline{->}(0,0)(30,0)
\psline{->}(0,0)(-20,30)

\rput[l](30,0){$2c-d$}
\rput[l](-23,35){$c$}
\rput[r](3,35){$b$}
\rput[l](-28,-35){$b+3c-3d$}
\rput[r](23,23){$d$}
\rput[l](-8,-1){$a$}
\end{pspicture}
\begin{pspicture}(-35,-50)(100,50)

\psline{->}(0,0)(0,30)
\psline{->}(0,0)(0,-30)
\psline{->}(0,0)(30,0)
\psline{->}(0,0)(-20,30)

\rput[l](30,0){$2z-w$}
\rput[l](-23,35){$z$}
\rput[r](3,35){$y$}
\rput[l](-28,-35){$y+3z-3w$}
\rput[r](23,23){$w$}
\rput[l](-8,-1){$x$}
\end{pspicture}
\]
The first weight has value $d$ on all maximal cones, and the second has value $w$. To calculate the product, we need a Riemann-Roch matrix. As detailed in Section \textcolor{black}{3}, we require a complete flag. In two dimensions, this is merely the data of a vector, so we pick the vector $(1,1)$. The resulting Riemann-Roch matrix is 
\[
\begin{pmatrix}
1     & 0 & 0 & 0 & 0 & 0 & 0 & 0 & 0 \\
\frac{1}{2}  & 1 & 0 & 0 & 0 & 0 & 0 & 0 & 0 \\
\frac{1}{2}  & 0 & 1 & 0 & 0 & 0 & 0 & 0 & 0 \\
\frac{1}{2}  & 0 & 0 & 1 & 0 & 0 & 0 & 0 & 0 \\
\frac{1}{2}  & 0 & 0 & 0 & 1 & 0 & 0 & 0 & 0 \\
\frac{5}{12} & \frac{1}{2} & \frac{1}{2} & 0 & 0 & 1 & 0 & 0 & 0 \\
\frac{1}{12} & 0 & \frac{1}{2} & \frac{1}{2} & 0 & 0 & 1 & 0 & 0 \\
\frac{41}{60} & 0 & 0 & \frac{1}{2} & \frac{1}{2} & 0 & 0 & 1 & 0 \\
-\frac{11}{60} & \frac{1}{2} & 0 & 0 & \frac{1}{2} & 0 & 0 & 0 & 1 \\
\end{pmatrix}.
\]
We must also select a displacement vector $v$ which specifies the values of $m^\alpha_{\beta,\gamma}$. If we choose e.g. $v=(5,1)$, then the non-zero $m^\alpha_{\beta,\gamma}$ are $m^{0}_{\rho_1,\rho_4}=2$, and $m^{0}_{\rho_2,\rho_3}=m^{0}_{0,\sigma_{34}}=m^{0}_{\sigma_{12},0}=m^{\rho_1}_{\rho_1,\sigma_{14}}=m^{\rho_1}_{\sigma_{12},\rho_1}=m^{\rho_2}_{\rho_2,\sigma_{23}}=m^{\rho_2}_{\sigma_{12},\rho_2}=m^{\rho_3}_{\rho_3,\sigma_{34}}=m^{\rho_3}_{\sigma_{23},\rho_3}=m^{\rho_4}_{\rho_4,\sigma_{34}}=m^{\rho_4}_{\sigma_{14},\rho_4}=1$.
The resulting weight is
\[
\begin{pspicture}(-55,-60)(100,70)

\psline{->}(0,0)(0,60)
\psline{->}(0,0)(0,-60)
\psline{->}(0,0)(60,0)
\psline{->}(0,0)(-40,60)

\rput[l](60,0){\footnotesize{$2cw+2dz-3dw$}}
\rput[l](-83,65){\footnotesize{$cw+dz-dw$}}
\rput[l](-5,65){\footnotesize{$bw+dy-dw$}}
\rput[l](-48,-65){\footnotesize{$bw+3(cw-dz)-7dw$}}
\rput[r](43,43){\footnotesize{$dw$}}
\rput[l](-8,-1){$ $}
\end{pspicture}
\]
\\
with a value of $a w - 2 b w - 8 c w + 9 d w + d x + 2 c y - 2 d y + 
 2bz + 6 cz - 8 dz$ on the origin.
\end{ex}

\section{Maps to Grothendieck weights}\label{sec:mapstogw}

The most straightforward relationship to describe is the map from Minkowski weights to Grothendieck weights:

\begin{prop}
Let $\Delta$ be a complete fan. There is an induced map $T:\MW^*(\Delta)\rightarrow\GW(\Delta)_\QQ$. For a simplicial fan this sends a weight $f\in\MW^k(\Delta)$ to the function $g$, defined by 
\[
g(\alpha) = \sum_{\alpha\prec\beta}\mu_\alpha(\beta) f(\beta).
\]
\end{prop}

This is immediate from the description of the Riemann-Roch transformation $\tau_X:K_\circ(X)\rightarrow A_*(X)_\QQ$ given in Proposition \ref{prop:rrmtx}, and the isomorphisms between $\MW^*(\Delta)$ and $A_*(X)^\vee$ when $\Delta$ is complete. For arbitrary toric varieties, one has a map to $\GW(\Delta)$ from functions on $\Delta$ obtained from elements of $A_*(X)^\vee$.

We compile a few comments about the relationship between Grothendieck weights and Minkowski weights in the following remark:
\begin{rem}
It is possible to use this proposition to algorithmically calculate the inverse image under $T$ of a Grothendieck weight. Let $g\in\GW(\Delta)$ be a Grothendieck weight, and suppose that $g\in I_k$. Then, the function on cones of dimension $k$ obtained by restriction of $g$ is a Minkowski weight: if $k$ is not the largest number that this is true, then $g|_k$ is uniformly $0$, and if $k$ is the largest, then this can be seen by Theorem \ref{thm:balancing}, using the fact that $\nu_\beta(\beta)=1$. Then, $g-T(g|_k)$ is an element of $\GW(\Delta)_\QQ$, and is in $I_{k+1}$. Inductively one obtains $T^{-1}(g)$. 

It is easy to check that the Riemann-Roch map sends $f\in\MW^k(\Delta)$ to $T(f)=g\in I_k\otimes\QQ$ satisfying $g|_{\Delta^{(k)}}=f$. It is natural to wonder whether, given any Minkowski weight $f\in\MW^k(\Delta)$, there is some $g\in I_k$ that satisfies the same equality when restricted to $\Delta^{(k)}$. This seems interesting because of the following proposition.
\end{rem}
For $f\in\MW^k(\Delta)$ we say that $g\in\GW(\Delta)$ lifts $f$ if $g(\alpha)=0$ for $\codim(\alpha)<k$ and $g|_{\Delta^{(k)}}=f$. Let $F_i$ be the $i$-th piece of the dimension filtration on the Grothendieck group, meaning that it is generated by coherent sheaves with support of dimension at most $i$. We also find it convenient to extend the definition of Minkowski weights for the next proposition. If $\Delta$ is not complete, one may still make sense of $\MW^k(\Delta)$ by defining it to be the set of functions on $\Delta^{(k)}$ given by $\alpha\rightarrow f([V(\alpha)])$ for $f$ a linear form on $A_k(X)$.
\begin{prop}
Let $\Delta$ be an arbitrary fan. Suppose that there exists $f\in\MW^k(\Delta)$ with no lift in $\GW(\Delta)$. Then $F_{k}$ is not saturated as a subgroup of $K_\circ(X)$.
\end{prop}
\noindent {\sl Proof:} 
We have the exact sequence 
\[
0\rightarrow F_{k}/F_{k-1}\rightarrow K_\circ(X)/F_{k-1}\rightarrow K_\circ(X)/F_{k}\rightarrow 0.
\]
The long exact sequence obtained after applying $(-)^\vee = Hom_\ZZ(-,\ZZ)$ is:
\[
\begin{tikzcd}
\scalemath{0.8}{0} \arrow[r]	& \scalemath{0.8}{(K_\circ(X)/F_{k})^\vee} \arrow[r] & \scalemath{0.8}{(K_\circ(X)/F_{k-1})^\vee} \arrow[r] & \scalemath{0.8}{(F_{k}/F_{k-1})^\vee} \ar[overlay,out=0, in=180]{dll}\\
 		&\scalemath{0.8}{\Ext^1_\ZZ(K_\circ(X)/F_{k},\ZZ)} \arrow[r] & \scalemath{0.8}{\Ext^1_\ZZ(K_\circ(X)/F_{k-1},\ZZ)} \arrow[r] & \scalemath{0.8}{\Ext^1_\ZZ(F_{k}/F_{k-1},\ZZ)}.
\end{tikzcd}
\]
Let us consider the first few terms. $(K_\circ(X)/F_{k-1})^\vee$ may be naturally identified with the ideal $I_{k}$ of Grothendieck weights which vanish on cones of codimension less than $k$, as defined in Corollary \ref{cor:filtration}. On the other hand, $F_{k}/F_{k-1}$ is the $k$-th piece of graded $K$-theory, and the map $A_{k}(X)\rightarrow F_{k}/F_{k-1}$ sending $[V]$ to $[\mathcal{O}_V]$ is an isomorphism after tensoring with $\QQ$ (see \cite[Chapter 18]{F98}). Thus, $(F_{k}/F_{k-1})^\vee\cong MW^{k}(\Delta)$. Thus the first few terms in the exact sequence become: 
\[
0 \rightarrow I_{k+1} \rightarrow I_{k} \rightarrow \MW^{k}(\Delta) \rightarrow\ldots .
\]
Suppose that Minkowski weight $f\in \MW^{k}(\Delta)$ does not have a lift. Then $I_{k}$ cannot surject onto $\MW^{k}(\Delta)$ otherwise the preimage of $f$ would be a lift. Thus the group $\Ext^1_\ZZ(K_\circ(X)/F_{k},\ZZ)$ cannot be trivial. But this group can (non-canonically) be identified with the torsion subgroup of $K_\circ(X)/F_{k}$. This being non-trivial is equivalent to $F_{k}$ not being saturated as a subgroup of $K_\circ(X)$.
\qed
\\

When $\Delta$ is complete, $\GW(\Delta)$ also admits a map from $\PExp(\Delta)$, which is the ring of continuous functions on $\Delta$ that are given on each cone $\alpha\in\Delta$ by an exponential function in the torus variables. Anderson and Payne showed that this ring is naturally isomorphic to $\opk^\circ_T(X)$, and thus the map to $\GW(\Delta)$ is induced by the forgetful map $\opk^\circ_T(X)\rightarrow\opk^\circ(X)$. We will require $K$-theoretic equivariant multiplicities $\epsilon_p^K (V(\alpha))$, where $p\in X^T$. These have been recently introduced in \cite{AGP}. They satisfy $\sum_{p\in X^T}\epsilon_p(V(\alpha))[i_{p*}(\mathcal{O}_p)]=[\mathcal{O}_{V(\alpha)}]$, where $i_{p*}$ is the pushforward in $K_\circ^T$ along the inclusion of the fixed point $p$.

\begin{thm}\label{thm:forgetful}
Let $\Delta$ be a complete fan. There is a commuting square
\[
\begin{tikzcd}
\opk^\circ_T(X(\Delta)) 	\arrow[r,"\cong"]\arrow[d,"\forgetful"]	&	\PExp(\Delta)\arrow[d,"\forgetful"]\\
\opk^\circ(X(\Delta))			\arrow[r,"\cong"]				& 	\GW(\Delta)
\end{tikzcd}
\]
where the forgetful map from $\PExp(\Delta)$ to $\GW(\Delta)$ sends a piecewise-exponential function $\phi$ to the limit of the function
\[
\alpha \rightarrow \sum_{\sigma\in\Delta^{(0)}}\epsilon^K_{V(\sigma)}(V(\alpha))f|_{\sigma},
\]
as the torus parameters approach $0\in N$.
\end{thm}
\noindent {\sl Proof:} 
If $\phi$ is our piecewise exponential function, then the corresponding $R(T)$-linear function $\psi:K^T_\circ(X)\rightarrow R(T)$ can be written explicitly via the projection formula:
\begin{align*}
\psi([\mathcal{O}_{V(\alpha)}])	& = \psi(\sum_{p\in X^T}\epsilon^K_p(V(\alpha))[i_{p*}(\mathcal{O}_p)]),\\
						& = i^*_{X^T}\psi(\sum_{p\in X^T}\epsilon^K_p(V(\alpha))[\mathcal{O}_p]),\\
						& = \sum_{p\in X^T}\epsilon^K_p(V(\alpha)) i^*_p\psi([\mathcal{O}_p]) = \sum_{p\in X^T}\epsilon^K_p(V(\alpha)) \phi|_{\sigma_p},
\end{align*}
where $\sigma_p$ is the maximal cone corresponding to $p$. 

Then, the forgetful map from $\opk^\circ_T(X)$ to $\opk^\circ(X)$ is induced by the projection $X\times T\rightarrow X$, meaning it is the pullback from $\opk^\circ_T(X)$ to $\opk^\circ_T(X\times T)\cong\opk^\circ(X)$. Via the identification of $\opk^\circ_T(X)$ with $R(T)$-linear maps from $K_\circ^T(X)$ to $R(T)$, and $\opk^\circ(X)$ with $K_\circ(X)^\vee$, the forgetful map sends $\psi:K_\circ^T(X)\rightarrow R(T)$ to the linear function on $K_\circ(X)$ sending $[\mathcal{O}_{V(\alpha)}]$ to the equivalence class in $\ZZ$ of $\psi([\mathcal{O}_{V(\alpha)}])$ (see the appendix of \cite{AGP} for more details). This is the same as taking the limit as the torus parameters approach $0\in N$.
\qed

\begin{ex}\label{ex:nonsurjsurface}
We apply this theorem to the toric variety $X$ with fan $\Delta$ in $N=\ZZ^2$ generated $(\pm 1,\pm 1)$. This example and the following corollary are analogous to \cite[Example 4.1 and Theorem 1.5]{KP}. In this case, a generating set for $\PExp(\Delta)$ over $R(T)$ for is given by the following functions:

\hspace{30 mm}
\begin{pspicture}(-35,-50)(100,50)

\psline{->}(0,0)(30,30)
\psline{->}(0,0)(-30,-30)
\psline{->}(0,0)(-30,30)
\psline{->}(0,0)(30,-30)

\rput[l](30,0){$1$}
\rput[r](-30,0){$1$}
\rput[l](-3,30){$1$}
\rput[l](-3,-30){$1$}

\end{pspicture}
\begin{pspicture}(-35,-50)(100,50)

\psline{->}(0,0)(30,30)
\psline{->}(0,0)(-30,-30)
\psline{->}(0,0)(-30,30)
\psline{->}(0,0)(30,-30)

\rput[l](30,0){$1-e^{x-y}$}
\rput[r](-30,0){$1-e^{x+y}$}
\rput[l](-3,30){$0$}
\rput[l](-14,-30){$1-e^{2x}$}

\end{pspicture}

\hspace{30 mm}
\begin{pspicture}(-35,-50)(100,50)

\psline{->}(0,0)(30,30)
\psline{->}(0,0)(-30,-30)
\psline{->}(0,0)(-30,30)
\psline{->}(0,0)(30,-30)

\rput[l](30,0){$0$}
\rput[r](-15,0){$(1-e^{x+y})(e^x-e^y)$}
\rput[l](-3,30){$0$}
\rput[l](-3,-30){$0$}

\end{pspicture}
\begin{pspicture}(-35,-50)(100,50)

\psline{->}(0,0)(30,30)
\psline{->}(0,0)(-30,-30)
\psline{->}(0,0)(-30,30)
\psline{->}(0,0)(30,-30)

\rput[l](30,0){$0$}
\rput[r](-30,0){$1-e^{x+y}$}
\rput[l](-3,30){$0$}
\rput[l](-19,-30){$1-e^{x+y}$}

\end{pspicture}

The images of these weights are the following Grothendieck weights.

\hspace{30 mm}
\begin{pspicture}(-35,-50)(100,50)

\psline{->}(0,0)(30,30)
\psline{->}(0,0)(-30,-30)
\psline{->}(0,0)(-30,30)
\psline{->}(0,0)(30,-30)

\rput[l](30,0){$1$}
\rput[r](-30,0){$1$}
\rput[l](-3,30){$1$}
\rput[l](-3,8){$1$}
\rput[l](-3,-30){$1$}
\rput[l](32,-32){$1$}
\rput[r](-32,-32){$1$}
\rput[r](-32,32){$1$}
\rput[l](32,32){$1$}

\end{pspicture}
\begin{pspicture}(-35,-50)(100,50)

\psline{->}(0,0)(30,30)
\psline{->}(0,0)(-30,-30)
\psline{->}(0,0)(-30,30)
\psline{->}(0,0)(30,-30)

\rput[l](30,0){$0$}
\rput[r](-30,0){$0$}
\rput[l](-3,30){$0$}
\rput[l](-3,-30){$0$}

\rput[l](-3,8){$2$}

\rput[l](32,-32){$-1$}
\rput[r](-32,-32){$1$}
\rput[r](-32,32){$-1$}
\rput[l](32,32){$1$}

\end{pspicture}

\hspace{30 mm}
\begin{pspicture}(-35,-50)(100,50)

\psline{->}(0,0)(30,30)
\psline{->}(0,0)(-30,-30)
\psline{->}(0,0)(-30,30)
\psline{->}(0,0)(30,-30)

\rput[l](30,0){$0$}
\rput[r](-30,0){$0$}
\rput[l](-3,30){$0$}
\rput[l](-3,-30){$0$}

\rput[l](-3,8){$2$}

\rput[l](32,-32){$0$}
\rput[r](-32,-32){$0$}
\rput[r](-32,32){$0$}
\rput[l](32,32){$0$}

\end{pspicture}
\begin{pspicture}(-35,-50)(100,50)

\psline{->}(0,0)(30,30)
\psline{->}(0,0)(-30,-30)
\psline{->}(0,0)(-30,30)
\psline{->}(0,0)(30,-30)

\rput[l](30,0){$0$}
\rput[r](-30,0){$0$}
\rput[l](-3,30){$0$}
\rput[l](-3,-30){$0$}

\rput[l](-8,10){$-1$}

\rput[l](32,-32){$1$}
\rput[r](-32,-32){$0$}
\rput[r](-32,32){$1$}
\rput[l](32,32){$0$}

\end{pspicture}

We go through the calculation for the piecewise exponential function at the top right: since the equivariant multiplicity of a point is just $1$, the value of the Grothendieck weight on any maximal cone is just the value of the piecewise exponential function at $0$. For the function we are considering this is $0$. For the ray $\rho$ generated by $(1,1)$, $V(\rho)$ is a $\PP^1$, and at the fixed point corresponding to the maximal cone $\sigma$ generated by $(1,1)$ and $(1,-1)$ the character on the tangent space is $y-x$, so the equivariant multiplicity is $\frac{1}{1-e^{x-y}}$, by \cite[Proposition 6.3]{AGP}. At the other fixed point of $V(\rho)$, the character is $x-y$, and so the multiplicity is $\frac{1}{1-e^{y-x}}$. The value of the Grothendieck weight on $\rho$ is then the limit of $\frac{0}{1-e^{y-x}}+ \frac{1-e^{x-y}}{1-e^{x-y}}$ as $x$ and $y$ approach $0$, which is $1$. Similarly, one gets that the value of the Grothendieck weight on the ray $(1,-1)$ is $-1$. The balancing conditions for Grothendieck weights determine the values on the other rays.

For the cone $\{0\}$, we will require the equivariant multiplicities $\epsilon_p(X)$. Since $X$ is singular at each fixed point, we can compute the equivariant multiplicity at the fixed point $p$ corresponding to $\sigma$ by resolving, e.g. by adding the ray $(1,0)$, and then summing over the new fixed points which map to $p$. One gets
\[
\epsilon_p(X)=\frac{1}{(1-e^y)(1-e^{x-y})}+\frac{1}{(1-e^{-y})(1-e^{x+y})}=\frac{1+e^x}{(1-e^{x+y})(1-e^{x-y})}.
\]
Let the ``bottom" fixed point be $q$, and the ``left" one be $r$. Then 
\begin{align*}
\epsilon_q(X) & =\frac{1+e^{-y}}{(1-e^{x-y})(1-e^{-x-y})},\\
\epsilon_r(X) & =\frac{1+e^{-x}}{(1-e^{-x-y})(1-e^{y-x})}.
\end{align*}
By the last theorem, $\alpha$ must be sent to the limit of 
\[
\frac{(1-e^{x-y})(1+e^x)}{(1-e^{x+y}))(1-e^{x-y})}+\frac{(1-e^{2x})(1+e^{-y})}{(1-e^{x-y})(1-e^{-x-y})}+\frac{(1-e^{x+y})(1+e^{-x})}{(1-e^{-x-y})(1-e^{y-x})},
\] as the torus parameters approach $0$, which is $2$.
\end{ex}
In fact, this example shows:
\begin{cor}\label{cor:nonsurj}
There exists a complete toric surface with a vector bundle with no finite length resolution by $T$-equivariant vector bundles.
\end{cor}

\noindent {\sl Proof:} 
In Example \ref{ex:nonsurjsurface}, the $\ZZ$-linear span of the Grothendieck weights calculated does not include the Grothendieck weight with $1$ at the origin and $0$ elsewhere, so $\PExp(\Delta)$ does not surject onto $\GW(\Delta)$. Thus, the forgetful map from $\opk^\circ_T(X)$ to $\opk^\circ(X)$ is not surjective. Since vector bundles induce linear forms on coherent sheaves by tensor product followed by pushforward to a point, there is a commutative square:
\[
\begin{tikzcd}
K_T^\circ(X)\arrow[r]\arrow[d] & \opk^\circ_T(X)\arrow[d]\\
K^\circ(X)\arrow[r] & \opk^\circ(X)
\end{tikzcd}.
\]
We know from \cite[Proposition 7.4]{AP} that the bottom map is surjective. Comparing the two ways of traversing the diagram, one sees that the map $K^\circ_T(X)\rightarrow K^\circ(X)$ cannot be surjective. This proves the corollary.
\qed

\section{Sch\" on subvarieties}
We recall some definitions from \cite{Te}: Let $Y$ be a closed subvariety of $T:=(\CC^*)^n$, and suppose that $X$ is a toric variety partially compactifying $T$, with fan $\Delta$. Let $\overline{Y}$ be the closure of $Y$ in $X$. The variety $\overline{Y}$ is called a \textit{tropical compactification} of $Y$ if the natural map $\Psi:T\times\overline{Y}\rightarrow X$ is faithfully flat and $\overline{Y}$ is complete. We call $X$ the \textit{ambient toric variety}. If $\Psi$ is also smooth, then $Y$ is called \textit{sch\" on}. In this case, \cite[Theorem 1.4]{Te} ensures us that $\overline{Y}$ is regularly embedded in $X$. In this case, the sum $\sum_{i=0}(-1)^i [\Tor_{\mathcal{O}_{X}}^i(\mathcal{O}_{\overline{Y}},\mathscr{F})]$ is well-defined for $\mathscr{F}$ any coherent sheaf, and so defines a element $[\mathcal{O}_{\overline{Y}}]\cap-$ in $K_\circ(X)^\vee$.

\begin{Defi}
Let $\Delta$ be the fan corresponding to $X$. Let $g_{\overline{Y}}:\Delta\rightarrow\ZZ$ be the Grothendieck weight associated to $[\mathcal{O}_{\overline{Y}}]\cap-$. We call $g_{\overline{Y}}$ the \textit{associated ($K$-theoretic) cocycle}.
\end{Defi}

Different tropical compactifications of the same sch\" on subvariety induce classes in different ambient toric varieties, which are compatible. To make this precise, we describe some functorial maps of Grothendieck weights. We denote the relative interior of a cone $\alpha$ by $\relint(\alpha)$:

\begin{lem}\label{lem:stability}
If $\phi:U\rightarrow X$ is an open $T$-equivariant inclusion of toric varieties, then there is an induced pushforward on Grothendieck weights, which is extension by $0$. If $\phi:X(\Delta^{'})\rightarrow X(\Delta)$ is a proper toric birational morphism, then there is an induced pullback of Grothendieck weights which sends $g\in\GW(\Delta)$ to the map $\Delta^{'}\rightarrow\ZZ$ sending $\alpha^{'}\in\Delta^{'}$ to $g(\alpha)$ for any $\alpha\in\Delta$ such that $\relint(\alpha)\cap \relint(\alpha')\neq \emptyset$.
\end{lem}

\noindent {\sl Proof:} 
The first case follows immediately from the fact that the natural induced map $\phi^*:K_\circ(X)\rightarrow K_\circ(U)$ is simply given by restriction. In the second case, let $\Delta'$ be a subdivision of $\Delta$ inducing a map $\phi:X'\rightarrow X$ of toric varieties. Suppose for $\alpha'\in\Delta', \alpha\in\Delta$, $\relint(\alpha')\cap\relint(\alpha)\neq\emptyset$. Then standard vanishing results for toric varieties imply that $\phi_*([\mathcal{O}_{V(\alpha')}])=[\mathcal{O}_{V(\alpha)}]$, so if we map $e_{\alpha'}\in\QQ^{\Delta'}$ to $e_{\alpha}\in\QQ^\Delta$, we have a map of exact sequences:
\[
\begin{tikzcd}
0\arrow[r] & Rel_{K_0(X')} \arrow[r]\arrow[d] & \QQ^{\Delta'} \arrow[r]\arrow[d] & K_\circ(X') \arrow[r]\arrow[d] & 0\\
0\arrow[r] & Rel_{K_0(X)} \arrow[r] & \QQ^{\Delta} \arrow[r] & K_\circ(X) \arrow[r] & 0.
\end{tikzcd}
\]
All we need to do is show that the kernel of the middle map surjects onto the kernel of the last map. By the snake lemma, it is enough to show that $Rel_{K_\circ(X')}$ surjects onto $Rel_{K_\circ(X)}$. We denote the projection map from $\QQ^\Delta$ to $\QQ^{\Delta^{(l)}}$ by $\pi_{l}$. Define a filtration on $Rel_{K_\circ(X)}$ by $F_k Rel_{K_\circ(X)}=\bigcap_{l>k}\ker \pi_l$. Then, for example, $F_n Rel_{K_\circ(X)}$ is just $Rel_{K_\circ(X)}$, and $F_0 Rel_{K_\circ(X)}$ are elements of the form $e_{\sigma_1}-e_{\sigma_2}$ for $\sigma_1,\sigma_2$ maximal cones in $\Delta$.

We show by induction that $F_k Rel_{K_{\circ}(X')}$ surjects onto $F_k Rel_{K_\circ(X)}$. For $k=0$ this is clear, so we assume it is true for some $k_0$. Let 
\[
s=\left(\sum_{\alpha\in\Delta^{(k_0+1)}} a_\alpha e_\alpha\right)+\left(\sum_{\substack{\alpha\in\Delta^{(l)} \\ l\leq k_0}}a_\alpha e_\alpha\right),
\]
be in $F_{k_0+1}Rel_{K_\circ(X)}$. Then by the Riemann-Roch theorem we know that in the Chow group $\sum_{\alpha\in\Delta^{(k_0+1)}} a_\alpha [V(\alpha)]=0$. Since $\phi$ is an envelope by \cite[Lemma 1]{P}, the analogous map from $Rel_{A_k(X')}$ to $Rel_{A_k(X)}$ is surjective for any $k$, so we can find an expression $0=\sum_{\alpha'\in\Delta'^{(k_0+1)}}a_{\alpha'}[V(\alpha')]$ such the sum of $a_{\alpha'}$ for all $\alpha'$ subdividing $\alpha$ is $a_\alpha$. Then $0=\sum_{\alpha'\in\Delta'^{(k_0+1)}}a_{\alpha'}\tau_{X'}^{-1}([V(\alpha')])$ so there is a relation $r=\sum_{\alpha'\in\Delta'^{(k_0+1)}}a_{\alpha'}e_{\alpha'}+\ldots$ which by the Riemann-Roch theorem is in $F_{k_0+1}Rel_{K_\circ(X)}$ which does not quite map to $s$. However, its image and $s$ only differ by an element of $F_{k_0}Rel_{K_\circ(X)}$. By our induction hypothesis we are done. \qed

\begin{cor}
Let $Y$ be a sch\" on subvariety of $T$ and $Y_1$ and $Y_2$ two tropical compactifications  in toric varieties $X_1$ and $X_2$. Consider the function $\wt{g}_{Y_i}:N_\QQ\rightarrow\ZZ$ induced from $g_{Y_i}$ by sending $x\in \relint(\alpha)$ to $g_{Y_i}(\alpha)$, and sending $x$ not in any cone in $\Delta$ to $0$. Then $\wt{g}_{Y_1}=\wt{g}_{Y_2}$.
\end{cor}

\noindent {\sl Proof:} 
Given any two tropical compactifications $Y_1$ and $Y_2$ in toric varieties $X_1$ and $X_2$, one can find a third $Y_3$ by taking the closure of $Y$ in $X_3$ the toric variety induced by common refinement of the overlap of the fans of $X_1$ and $X_2$. The resulting toric variety will a proper birational morphism onto $T$-invariant open subsets of $X_1$ and $X_2$, which is an isomorphism on $T$. Additionally the fan of $X_3$ will contain $\Trop(Y)$, and so $Y_3$ will be a tropical compactification, whose associated cocycle can be related to the cocycles of $Y_1$ and $Y_2$ by maps as described in the previous proposition.
\qed

\begin{ex}
The moduli space of genus $0$ curves $M_{0,n}$ can be identified with $(\PP\smallsetminus\{0,1,\infty\})^{n-3}$. This can be realized as a sch\" on subvariety of $(\CC^*)^{\binom{n-1}{2}-1}$, by \cite{Te}. We recall the explicit description of the tropicalization of $M_{0,n}$ given in \cite{GKM}.

The abstract moduli space of tropical rational curves $\mathcal{M}_{0,n}$ is the set consisting of all connected trivalent metric trees with $n$ leaves $p_1,\ldots,p_n$. Let $Q_n$ be the quotient of $\RR^{\binom{n}{2}}$ by the image of $\RR^n$ under the map sending $(a_1,\ldots,a_n)$ to $(a_i+a_j)_{i,j}$. Then, one has a map from $\mathcal{M}_{0,n}$ to $Q^n$ which sends $C$ to $(dist_C(p_i,p_j))_{i,j}$, where $dist_C$ denotes the distance function on $C$. This map is in fact injective, and its image is the support of a rational polyhedral fan. This can be identified with the tropicalization of $M_{0,n}$, and the associated (non-compact) toric variety $X$ provides a tropical compactification of $M_{0,n}$, which can be identified with the Mumford-Knudsen compactification $\overline{M}_{0,n}$. Since the restriction of $\overline{M}_{0,n}$ to each $T$-invariant subvariety of $X$ is a rational variety, the associated cocycle is just the constant function $1$ on all the cones of $\Trop(M_{0,n})$.
\end{ex}

	
	
		%
	
We can use the associated cocycle of $\overline{Y}$ to calculate some of its invariants. Suppose that the ambient toric variety $X$ is complete and $E$ is a vector bundle on $X$. 
\begin{Defi}
Let the function $g_E:\Delta\rightarrow\ZZ$ which sends $\alpha$ to $\chi([E|_{V(\alpha)}])$ be called the associated cocycle of $E$.
\end{Defi}

Note that the function $g_E$ is a Grothendieck weight by definition. Then the following proposition is formal, but it allows us to reduce computation of Euler characteristics on a variety which is not necessarily toric, i.e. $\overline{Y}$, to computations of Euler characteristics on toric varieties.

\begin{prop}
If $E$ is a vector bundle on a complete ambient space $X$, the value of $g_{\overline{Y}}\cup g_{E}$ on the cone $\{0\}$ is equal to $\chi(E|_{\overline{Y}})$.
\end{prop}

\noindent {\sl Proof:} 
Since the ambient toric variety $X$ is complete, the rings $\GW(\Delta), K_\circ(X)^\vee$, and $\opk^\circ(X)$ are all isomorphic. Along this identification, $g_{\overline{Y}}$ corresponds to the $\opk^\circ(X)$ class $([-\otimes_{\mathcal{O}_{X^{'}}}^L Lf^{*}\mathcal{O}_{\overline{Y}}])_{f:X^{'}\rightarrow X}$, and $g_E$ corresponds to $([-\otimes_{\mathcal{O}_{X^{'}}} f^*E])_{f:X^{'}\rightarrow X}$. The product in $\opk^\circ(X)$ is merely composition of endomorphisms for each $f:X^{'}\rightarrow X$, so the composition of these classes is $([-\otimes_{\mathcal{O}_{X^{'}}}^L Lf^{*}\mathcal{O}_{\overline{Y}}\otimes_{\mathcal{O}_{X^{'}}} f^*E])_{f:X^{'}\rightarrow X}$. The corresponding $K_\circ(X)^\vee$ element is then $\chi(-\otimes^L_{\mathcal{O}_X}\mathcal{O}_{\overline{Y}}\otimes_{\mathcal{O}_X}E)$. To get the value of the Grothendieck weight at $\{0\}$, we simply evaluate this linear form on the structure sheaf of $V(\{0\})=X$, i.e. the class $[\mathcal{O}_X]$. But $\chi(\mathcal{O}_{\overline{Y}}\otimes_{\mathcal{O}_X}E)=\chi(E|_{\overline{Y}})$.
\qed

Inspired by the fact that a piecewise linear function (tropical Cartier divisor) on a tropical variety gives rise to a tropical cycle, e.g. \cite{AR}, we come to the following proposition, with our notation as earlier in this section.
\begin{thm}
Let $\Delta$ be complete, and let $F$ be the forgetful map from $\PExp(\Delta)$ to $\GW(\Delta)$. Let $I$ be the ideal of piecewise exponential functions on $\Delta$ vanishing on $\Trop (Y)$. Then $F(-)\cup g_{\overline{Y}}$ defines a map from $\PExp(\Delta)/I$ to the ideal of functions in $\GW(\Delta)$ supported on $\Trop(Y)$.  
\end{thm}
\noindent {\sl Proof:} 
Let $g\in I$. By Lemma \ref{lem:stability}, we may assume that $\Delta$ is smooth by pulling everything back to a smooth refinement. The subfan $\Trop(Y)$ in $\Delta$ defines a $T$-invariant open set $U$ in $X$ which is also a toric variety. By assumption, $g$ restricts to $0$ on $\Trop(Y)$, so its corresponding class in $\opk^\circ_T(X)$ must restrict to $0$ when pulled back to $U$. Since $X$ is smooth, there is an isomorphism between $K^T_\circ(X)$ and $\opk_T^\circ(X)$. Thus, we may apply the localization exact sequence for $K^T_\circ$ which says that classes which pull back to $0$ on $U$ must be the pushforward of some class on $X\smallsetminus U$. But $g_{\overline{Y}}$ evaluates such classes to $0$ since it is supported on $\Trop(Y)$. This shows that $g$ must map to $0$, so $F(-)\cup g_{\overline{Y}}$ descends to a map from $\PExp(\Delta)/I$ to $\GW(\Delta)$. The image must land in the ideal of Grothendieck weights supported on $\Trop(Y)$ because $g_{\overline{Y}}$ is in this ideal.
\qed

\begin{ex}
On $\Trop(M_{0,n})$, there are certain $\ZZ$-valued functions $\Psi_k$ called tropical Psi classes. These can be viewed as Minkowski weights on any toric compactification of $X$, whose corresponding Chow classes restrict to classical Psi classes on $\overline{M}_{0,n}$. In \cite{KM}, the authors found equivariant Chow classes $f_k$ on compactifications $X_k$ of $X$ whose corresponding non-equivariant classes restrict to $\binom{n-1}{k}\cdot\Psi_k$. The compactification $X_k$ has maximal cones generated by $\binom{n-1}{2}-1$-element subsets of $v_{i,j}$ for $\{i,j\}\in [n]\smallsetminus\{k\}$. Gathmann and Markwig define $f_k$ as the piecewise polynomial function on the fan of $X_k$ determined by 
\begin{equation*}
    f_k(v_{i,j}) = \begin{cases}
	       	0					& i\text{ or } j=k\\
		1 					& \text{otherwise.}
           \end{cases}
\end{equation*} 
For a line bundle on a toric variety, the (equivariant) Chern class will be a piecewise linear function. One can obtain the equivariant operational $K$-theory class of the line bundle by exponentiating this piecewise linear function. Thus, one may obtain the associated cocycle of the cotangent line bundle $L_k$ on $\overline{M}_{0,n}$ by computing the image of $\exp({\frac{f_k}{\binom{n-1}{2}}})$ in $\GW(\Delta)$. 

For $n=5$, the associated cocycle $g_{L_k}$ is the Grothendieck weight sending a cone $\beta$ to $\codim(\beta)+1$. By implementing Theorem \ref{thm:prod}, we calculated using Sage that $(g_{L_k}\cup g_{\overline{M}_{0,n}})(\{0\})=\chi(L_k,\overline{M}_{0,n})=3$.
\end{ex}

\subsection*{Acknowledgements:}
This research was partially supported by the NSF-RTG grant, \# DMS-1547357. I would like thank Ang\' elica Cueto and Eric Katz for comments on an embryonic version of the first few sections of this document, and Kiumars Kaveh and Sam Payne for comments on a more recent version. Most importantly, I thank David Anderson for the idea that one should consider a K-theoretic analogue of Minkowski weights in the first place, and for his comments and direction throughout the course of creating this document. 

\appendix

\section{Multiplicities of cones}
Given cones $\alpha\prec\beta$, we denote the \textit{relative multiplicity} of $\overline{\beta}$ in $N_\alpha$ by $\mult_\alpha(\beta)$. If $\alpha=\{0\}$ we obtain that $\mult_\alpha(\beta)=\mult(\beta)$ is the usual multiplicity of $\beta$. The following lemma describes relative multiplicities of simplicial cones in terms of usual multiplicities. Let $\alpha$ have rays $\rho_1,\ldots,\rho_k$, and $\beta$ have rays $\rho_1,\ldots,\rho_k,\rho_{k+1},\ldots,\rho_l$:

\begin{lem}
$\mult_\alpha(\beta) = \mult(\beta)\prod_{i=k+1}^l\frac{\mult(\alpha)}{\mult(\alpha+\rho_i)}$.
\end{lem}
\noindent {\sl Proof:} 
To simplify notation, we assume that $\beta$ is a maximal cone. Then we have the following diagram of exact sequences:
\[
\begin{tikzcd}
0 \arrow[r]& \langle v_{\rho_1},\ldots,v_{\rho_l} \rangle \arrow[r]\arrow[d] & N \arrow[d]\arrow[r]& N/\langle v_{\rho_1},\ldots,v_{\rho_l} \rangle \arrow[r]\arrow[d]& 0\\
0 \arrow[r]& \langle v_{\overline{\rho_{k+1}}},\ldots,v_{\overline{\rho_l}} \rangle \arrow[r] & N_\alpha \arrow[r]& N_\alpha/\langle v_{\overline{\rho_{k+1}}},\ldots,v_{\overline{\rho_l}} \rangle \arrow[r]& 0,
\end{tikzcd}
\]
where the top and bottom quotient groups on the right have cardinality $\mult(\beta)$ and $\mult_\alpha(\beta)$ respectively. We add the kernels and cokernels to the diagram:

\[
\begin{tikzcd}
0 \arrow[r]& \langle v_{\rho_1},\ldots,v_{\rho_k} \rangle \arrow[r,"\cong"]\arrow[d] & \langle v_{\rho_1},\ldots,v_{\rho_k} \rangle \arrow[d]\arrow[r]& A \arrow[d]\\
0 \arrow[r]& \langle v_{\rho_1},\ldots,v_{\rho_l} \rangle \arrow[r]\arrow[d] & N \arrow[d]\arrow[r]& N/\langle v_{\rho_1},\ldots,v_{\rho_l} \rangle \arrow[r]\arrow[d]& 0\\
0 \arrow[r]& \langle v_{\overline{\rho_{k+1}}},\ldots,v_{\overline{\rho_l}} \rangle \arrow[r]\arrow[d] & N_\alpha \arrow[r]\arrow[d]& N_\alpha/\langle v_{\overline{\rho_{k+1}}},\ldots,v_{\overline{\rho_l}} \rangle \arrow[r]\arrow[d]& 0\\
	& B \arrow[r] & 0 \arrow[r]& 0 \arrow[r]& 0.
\end{tikzcd}
\]
By the snake lemma, the sequence of kernels leading to cokernels is exact, so in fact $A\cong B$. Thus $\mult(\beta)=\mult_\alpha(\beta)|B|$. The cardinality of $B$ on the other hand is also easy to determine: the image of $\langle v_{\rho_1},\ldots,v_{\rho_l} \rangle$ in $\langle v_{\overline{\rho_{k+1}}},\ldots,v_{\overline{\rho_l}} \rangle$ is just $\langle \overline{v_{\rho_{k+1}}},\ldots,\overline{v_{\rho_l}}\rangle$. If $\overline{v_\rho}=b_\rho v_{\overline{\rho}}$, the cardinality of the cokernel (i.e. $B$) is $\prod_{i=k+1}^l b_{\rho_i}$. But in the proof of Proposition \ref{prop:multratio}, we saw $b_\rho=\frac{\mult(\alpha+\rho)}{\mult(\alpha)}$, which proves the claim.
\qed

\vskip 1cm

A.Shah.: 
Department of Mathematics,
Ohio State University,
Columbus, Ohio. \hfill\break
 Email:
{\tt     shah.1099@osu.edu  }

\ed